\documentclass[11pt,draft	]{article}
\usepackage{amsfonts}
\usepackage{amssymb}
\usepackage{amsmath}
\usepackage{epsfig,epic,eepic}

\usepackage{psfrag}

\newcommand{\R}{\mathbb R}
\newcommand{\N}{\mathbb N}

\newcommand{\hop}{\vskip .2cm\noindent}
\newcommand{\hip}{\vskip .1cm\noindent}

\newcommand{\hD}{\widehat D}
\newcommand{\ff}{\mathcal F}

\newcommand{\hM}{\widehat M}
\newcommand{\hg}{\widehat g}
\newcommand{\hX}{\widehat X}
\newcommand{\hY}{\widehat Y}
\newcommand{\hZ}{\widehat Z}
\newcommand{\wX}{X}
\newcommand{\wY}{Y}
\newcommand{\wZ}{Z}

\newtheorem{enonce}{}[section]

\newtheorem{thm}{Theorem}
\newtheorem{cor}[enonce]{Corollary}
\newtheorem{prop}[enonce]{Proposition}
\newtheorem{lem}[enonce]{Lemma}
\newtheorem{defi}[enonce]{Definition}
\newtheorem{fact}[enonce]{Fact}
\newtheorem{quest}[enonce]{Question}
\author {Pierre Mounoud}
\title{Gallot-Tanno theorem for closed incomplete pseudo-Riemannian manifolds and application.}
\date{}
\begin{document}
\maketitle
\begin{abstract}
In this article we extend the Gallot-Tanno theorem to closed pseudo-Riemannian manifolds. It is done by showing that if the cone over such a manifold  admits a parallel symmetric $2$-tensor then it is  incomplete and has non zero constant curvature. An application of this result  to the existence of metrics with distinct Levi-Civita connections but having the same unparametrized geodesics  is given.
\end{abstract}
\section{Introduction.}
Let $(M,g)$ be a pseudo-Riemannian manifold and let $D$ be its covariant derivative. Along this paper we will not consider Riemannian metrics (positive or negative) as pseudo-Riemannian one.
We are interested in the existence of non-constant solutions of the following equation:
$$(*)\qquad DDD \alpha (X,Y,Z)+ 2 (D\alpha  \otimes g)(X,Y,Z)+(D\alpha \otimes g)(Y,X,Z)+(D\alpha \otimes g)(Z,X,Y)=0,$$
where $\alpha$ is a function from $M$ to $\R$, and $X$, $Y$, $Z$ are vectors tangent to $M$.

This equation has already been studied, mostly in the Riemannian setting. For example Gallot and Tanno independently showed (see \cite{Ga} and \cite{Ta})  that if $(M,g)$ is a Riemannian complete manifold admitting a non-constant solution to equation $(*)$ then $(M,g)$ is a quotient of the round sphere. They were motivated by the fact that the eigenfunction related to the second eigenvalues of the Laplacian satisfies the equation $(*)$.

This equation also appears in the work of Solodovnikov \cite{So}, he showed that if a Riemannian metrics $g$ admits  lots of geodesically equivalent metrics some of them having a different Levi-Civita connection then there exists a real number $c$ and a non-constant solution of $(*)$ for the metric $cg$ (see section \ref{lichne} for definitions and more precise statements). This result has been recently extended to the pseudo-Riemannian setting by V. Kiosak and V.Matveev in \cite{Ma2}.

The main result of this article is the following generalization of the  Gallot-Tanno theorem:
\begin{thm}\label{pr}
 If $(M,g)$ is a closed (ie compact without boundary) pseudo-Riemannian (but not Riemannian) manifold which admits a non-constant solution to the equation 
$$(*)\qquad DDD \alpha (X,Y,Z)+ 2 (D\alpha  \otimes g)(X,Y,Z)+(D\alpha \otimes g)(Y,X,Z)+(D\alpha \otimes g)(Z,X,Y)=0,$$
then it is closed geodesically incomplete pseudo-Riemannian manifold of constant  curvature equal to $1$.
\end{thm}
This clearly implies (cf. theorem \ref{mob}) that if a closed pseudo-Riemannian metric of non-constant curvature admits lots of geodesically equivalent metrics then they share the same Levi-Civita connection.
Some elements from the statement of theorem \ref{pr} deserves to be commented.

Let us  start with the hypothesis of closeness --as the same result but under the additional assumption of geodesic completeness is given in \cite{Ma} and (part of it) in \cite{Leist} and as Gallot and Tanno assumed only completeness. In the realm of Riemannian geometry, the unit tangent bundle of a compact manifold being compact, closeness implies geodesic completeness. This is no more true in pseudo-Riemannian geometry, incomplete metrics on compact manifolds are abundant, for example one can deduce from the article \cite{CR} of Y. Carri\`ere and L. Rozoy that the set of incomplete Lorentzian $2$ dimensional tori is dense in the set of Lorentzian tori. It seems that completeness and closeness are quite independent properties. Moreover, example 3.1 of \cite {Leist} from Alekseevsky and al. consists of non compact complete pseudo-Riemannian manifolds admitting non-constant solutions to $(*)$.

We continue our comments with the  conclusion of the theorem. Contrarily to the Gallot-Tanno theorem we were not able to find a closed manifold on which the equation $(*)$ has a non trivial solution. The reason being that no closed incomplete pseudo-Riemannian manifold of constant  curvature equal to $1$ are known. It is  conjectured that they do not exist. 
Actually if we assume that the metric is Lorentzian (ie of signature $(n-1,1)$ or $(1,n-1)$) the Carri\`ere-Klingler theorem  tells us that they do not exist, therefore the conclusion of theorem \ref{pr} can be strengthened in that case (see corollary \ref{lorentz}): the equation $(*)$ does not have any solutions.

Finally, the Riemannian or Lorentzian oriented reader should not pay too much attention to the sign of the curvature. Indeed, as no assumption is made about the signature of $g$, the sign of the curvature is meaningless: if $g$ has constant curvature equal to $1$ then $-g$ has constant curvature~$-1$.

The organization of the article is as follows. In section \ref{section} the link between solutions of $(*)$ and parallel symmetric $2$-tensor on the cone over $(M,g)$  is given and is used to show that the existence of a solution of $(*)$ implies that the cone is decomposable. In section \ref{decompo} the decomposable case is studied and theorem \ref{pr} is proven.  At last, section \ref{lichne} is devoted to the application of theorem \ref{pr} to the question of the existence of metrics with distinct Levi-Civita connections but having the same unparametrized geodesics.
\section{Parallel symmetric $2$-tensors on the cone over a manifold.}\label{section}
We start by giving the definition of cones over a pseudo-Riemannian manifolds.
\begin{defi}
Let $(M,g)$ be a pseudo-Riemannian manifold. We call \emph{cone manifold} over $(M,g)$ the manifold $\hM=\R^*_+\times M$ endowed with the metric $\hg$ defined by $\hg=\,dr^2+ r^2 g$.
\end{defi}
We will denote by $D$ the Levi-Civita connection of $g$ and by $\hD$ the Levi-Civita connection of $\hg$.

 The holonomy of cones over pseudo-Riemannian is strongly related to the equation $(*)$ seen in the introduction. This relation is given by the following proposition, which is almost contained in the proofs of corollaire 3.3 from \cite{Ga} (for an implication)  and  corollary 1 in \cite{Ma} (for the reciprocal). As we will use  some lines from it, as those proofs have a non empty intersection,  and for the convenience of the reader, we give its proof but it does not pretend to be new.
\begin{prop}\label{GaMa}
Let $(M,g)$ be a pseudo-Riemannian manifold. 
Let  $(\hM,\hg)$ be the cone manifold over $(M,g)$.
There exists a smooth non-constant  function $\alpha: M\rightarrow \R$ such that for any vectorfields $X$, $Y$, $Z$ of $M$ we have:
$$(*)\qquad DDD \alpha (X,Y,Z)+ 2 (D\alpha  \otimes g)(X,Y,Z)+(D\alpha \otimes g)(Y,X,Z)+(D\alpha \otimes g)(Z,X,Y)=0,$$
if and only if there exists a non-trivial symmetric parallel $2$-tensor on $(\hM,\hg)$.

More precisely if $\alpha$ is a non-trivial solution of $(*)$ then the Hessian of the function $A:\hM\rightarrow \R$ defined by $A(r,m)=r^2\alpha(m)$ is parallel (ie $\hD\hD\hD A=0$). Conversely if $T$ is a symmetric parallel $2$-tensor on $\hM$ then $T(\partial_r,\partial_r)$ does not depend on $r$ and is a solution of $(*)$. Moreover $2T$ is the Hessian of the function $A_T$ defined by $A_T(r,m)=r^2T_{(r,m)}(\partial_r,\partial_r)$.
\end{prop}
{\bf Proof.}
Let $\alpha$ be a solution of $(*)$ and $A$ be as above. We follow the proof of \cite{Ga}, the fact that  $g$ is now pseudo-Riemannian does not affect it.

Let $m$ be a point of $M$ and $x$, $y$, $z$ be three vectors of $T_mM$. There exists  vector fields,  $X$, $Y$, $Z$ such that $X(m)=x$, $Y(m)=y$, $Z(m)=z$ and that $DX(m)=DY(m)=DZ(m)=0$. We will denote the same way their lift to $\hM$.
\begin{fact}\label{fact}
The Levi-Civita connection of $\hg$ is given by 
$$\hD_{\wX}\wY={D_XY}-rg(X,Y)\partial_r,\quad \hD_{\partial_r}\partial_r=0, \quad \hD_{\partial_r}\wX=\hD_{\wX}\partial_r=\frac{1}{r}\wX.$$
\end{fact}
{\bf Proof.} Using $\hg([\wX,\wY],\partial_r)=0$ and $[\partial_r,\wX]=[\partial_r,\wY]=0$, we have
$$2\hg(\hD_{\wX}\wY,\partial_r)=-\partial_r.\hg(\wX,\wY)=-2rg(X,Y).$$
Similarly we have that $2\hg(\hD_{\wX}\wY,\wZ)=r^2g(D_XY,Z)$. It implies the first assertion. 
The two others can be shown the same way.$\Box$
\hip
We define the vectorfields $\hX$, $\hY$, $\hZ$ by 
$$\hX=\frac{1}{r}\wX,\quad  \hY=\frac{1}{r}\wY,\quad \hZ=\frac{1}{r}\wZ.$$
From now on, even if it is not specified, the computations are made at a point above $m$. We deduce from  fact \ref{fact} that (at $m$)
$$\hD_{\partial_r}\partial_r=\hD_{\partial_r}\hX= \hD_{\partial_r}\hY=\hD_{\partial_r}\hZ=0,\qquad 
\hD_{\hX}\partial_r=\frac{1}{r}\hX\quad \mathrm{and} \quad \hD_{\hX}\hY=-\frac{1}{r}g(X,Y)\partial_r.$$
Hence, we have 
$$
\hD A(\hZ)=rD\alpha(Z) \quad \mathrm{and}\quad \hD A(\partial_r)=2r\alpha.$$
Then 
\begin{equation}\hD\hD A (\hZ,\partial_r)=\hD\hD A(\partial_r,\hZ)=\partial_r.(r\,D \alpha(Z))-\hD A(\hD_{\partial_r}\hZ)=D\alpha(Z), \label{deriv}\end{equation}
and similarly
\begin{equation}
\hD\hD A(\partial_r,\partial_r)=2\alpha \label{alpha}\end{equation}
Using that, $DY(m)=0$, we get $Y.D A(Z)=DD A(Y,Z)$ and
\begin{equation}\hD\hD A(\hY,\hZ)=\hY.(rD_Z\alpha)+\frac{1}{r}g(Y,Z)\hD A(\partial_r)=DD\alpha(Y,Z)+2g(Y,Z)\alpha.\label{hess}\end{equation}
The proof that $\hD^3 A(\partial_r,..)=\hD^3 A(.,\partial_r,.)=\hD^3 A(.,.,\partial_r)=0$ is left to the reader or can be found in \cite{Ga}. The last thing to check is
\begin{eqnarray}
 \hD^3 A(\wX,\hY,\hZ)&=&X.(\hD\hD A(\hY,\hZ))-\hD\hD(\hD_{\wX}\hY,\hZ)-\hD\hD(\hY,\hD_{\wX}\hZ)\nonumber \\ 
                        &=&X.(DD\alpha(Y,Z)+2g(Y,Z)\alpha)+\hD\hD A(\frac{1}{r}g(X,Y)\partial_r,\hZ)+\hD\hD A(\frac{1}{r}g(X,Z)\partial_r,\hY)\nonumber \\
                         &=&DDD \alpha(X,Y,Z)+2g(Y,Z) D\alpha(X)+g(X,Y)D\alpha(Z)+g(X,Z)D\alpha(Y).\label{last}
\end{eqnarray}
This prove the first half of the proposition. 

Let $T$ be a parallel $2$-tensor on $(\hM,\hg)$. We have
$$ 0=\hD T (\partial_r,\partial_r,\partial_r)=\partial_r.T(\partial_r,\partial_r)-2T(\hD_{\partial_r}\partial_r,\partial_r).$$
As $\hD_{\partial_r}\partial_r=0$, we have $\partial_r.T(\partial_r,\partial_r)=0$. Thus $T(\partial_r,\partial_r)$ is a function on $M$. Now we follow  \cite{Ma} to prove that this function is a solution of the equation $(*)$.
\begin{lem}[see \cite{Ma}]\label{retour}
 Let $T$ be a symmetric parallel $2$-tensor on $(\hM,hg)$, and let $\alpha$ be the function defined by $T(\partial_r,\partial_r)$. Let $X$, $Y$, $Z$ and $\hX$, $\hY$, $\hZ$ be has above. We have 
$$\begin{array}{rl}
   2T(\partial_r,\hX)=&D\alpha(X)\\
   2T(\hX,\hY)=&2g(X,Y)\alpha+DD\alpha(X,Y).
  \end{array}
$$
\end{lem}
{\bf Proof.} We start from $\hD T=0$, using fact \ref{fact} we have:
\begin{equation}
0=\hD T(\hX,\partial_r,\partial_r)=D\alpha(\hX)-2T(\hD_{\hX}\partial_r,\partial_r)=\frac{1}{r}\big(D\alpha(X)-2T(\hX,\partial_r)\big).\label{up}
\end{equation}
This shows the first assertion.
The second one is shown the same way, using again \ref{fact} and \ref{up}
$$\begin{array}{rl}
0=\hD T(\hX,\hY,\partial_r)=&\hX.T(Y,\partial_r)-T(\hD_{\hX}\hY,\partial_r)-T(\hY,\hD_{\hX}\partial_r)\\
   =&\frac{1}{2r}\big(DD\alpha(X,Y)+2g(X,Y)\alpha - 2T(\hY,\hX)\big).\qquad \Box
\end{array}$$

Comparing the result of lemma \ref{retour} and (\ref{deriv}), (\ref{alpha}), (\ref{hess}), we see that $2T$ is the Hessian of $A=r^2\alpha$. The relation $\hD T=0$ becomes $\hD^3A=0$ and the equation (\ref{last}) says that $\alpha$ is a solution of $(*)$. $\Box$
\begin{cor}\label{coco}
 If $\alpha$ is a solution of $(*)$ which is constant on an open subset $U$ of $M$ then $\alpha$ is constant on $M$.
\end{cor}
{\bf Proof.} As for any $k\in \R$, the function $\alpha+k$ is also a solution of $(*)$, we can assume that for any $m\in U$, $\alpha(m)=0$. As $D \alpha(m)=0$, it follows from  (\ref{deriv}) and (\ref{alpha}) that, for any $r>0$ and any $m\in U$, $\hD\hD A_{(r,m)} (\partial_r,.)$ vanishes on $TM$ and takes the value $2\alpha(m)$ on $\partial_r$. It means that
\begin{equation}\hD\hD A_{(r,m)} (\partial_r,.)=2\alpha(m)g_{(r,m)}(\partial_r,.)=0.\label{decadix}\end{equation}
Moreover as for any $m\in U$, we have $DD\alpha(m)=0$ then 
$$\hD\hD A_{(r,m)}(\hY,\hZ)=2g(Y,Z)\alpha(m)=0.$$
Hence the Hessian of $A$ vanishes on $\R^*_+\times U$, as it is parallel it vanishes everywhere. The gradient of $A$ is therefore parallel but it vanishes also on $\R^*_+\times U$ therefore $\alpha$ is constant.$\Box$
\hip
Contrarily to the Riemannian case, the existence of a parallel symmetric $2$-tensor on a pseudo-Rieman\-nian manifold does not implie that the manifold is decomposable (ie that it possess parallel non-degenerate distributions). It is a consequence of the fact that the self-adjoint endomorphism associated to such a tensor and the metric can not always be simultaneously diagonalized (see \cite{Bouba} for more precise results). However, the situation is more simple for cones as shows the following.
\begin{prop}\label{para->decomp}
 Let $(M,g)$ be a closed pseudo-Riemannian manifold. If the equation $(*)$ has a non-constant solution then $(\hM,\hg)$ is decomposable.
\end{prop}
{\bf Proof.}
Let $\alpha$ be a solution of $(*)$ on $M$. As $M$ is closed there exists two critical points $m_-$ and $m_+$ of $\alpha$  associated to distinct critical values (ie $d\alpha(m_{\pm})=0$ and $\alpha(m_-)\neq\alpha(m_+)$). Equation (\ref{decadix}) tells us that  we have:
$$\hD\hD A_{(r,m_{\pm})} (\partial_r,.)=2\alpha(m_{\pm})g_{(r,m_{\pm})}(\partial_r,.).$$ It means that the self-adjoint endomorphism  associated to $\hD\hD A$ has two distinct real eigenvalues.

 Furthermore it is well known that  the characteristic spaces of a self-adjoint endomorphism provide an orthogonal (and therefore non-degenerate) decomposition of the tangent bundle of $M$. In our case, the endomorphism being  moreover parallel, this decomposition is also parallel. The number of eigenvalues being greater than $1$, $(\hM,\hg)$ is therefore decomposable.$\Box$
\hip
Proposition \ref{para->decomp} does not say that a cone with interesting holonomy is automatically decomposable. For example, the cone may admit anti-symmetric parallel $2$-tensors. The reader can consult \cite{Leist} for a more systematic study of the holonomy of cones.
\section{Decomposable cones.}\label{decompo}
 Let $(M,g)$ be a \emph{closed} pseudo-Riemannian manifold such that  its cone $(\hM,\hg)$ admits  two non degenerate complementary parallel distributions $V_1$ and $V_2$. We recall that parallel distributions are integrable and we will denote by $\ff_i$ the foliation spanned by $V_i$.

Let $T_1$ and $T_2$ be the symmetric  $2$-tensors on $\hM$ defined for $i\in\{1,2\}$ by
$$T_i(u,v)=\hg(u_i,v_i),$$
where $u_i$ and $v_i$ are the factors of the decomposition of $u$ and $v$ according to the splitting $T\hM=V_1\oplus V_2$.
The distributions $V_i$ are parallel then the tensors $T_i$ are also parallel.

We set $\alpha=T_1(\partial_r,\partial_r)$, as $\hg(\partial_r,\partial_r)=1$ we have $1-\alpha=T_2(\partial_r,\partial_r)$.
As in section \ref{section} we define on $\hM$  the functions $A^1$ and $A^2$ by $A^1(r,m)=r^2\alpha(m)$ and $A^2(r,m)=r^2(1-\alpha(m))$.

From proposition \ref{GaMa}  and corollary \ref{coco} we deduce the folowing proposition (proven also in \cite{Leist}).
\begin{prop}
\begin{enumerate}
\item $\alpha$ is a function on $M$ and it is a solution of the equation $(*)$.
\item The open subset $U=\{m\in M\,|\, \alpha(m)\neq 0,1\}$ is dense.
\end{enumerate}
\end{prop}
This proposition has the following corollary.
\begin{cor}\label{critic}
For any $m\in M$, we have $0\leq \alpha(m)\leq 1$. In particular, any point $m$ such that $\alpha(m)=0$ or $1$ is a critical point. Moreover $0$ and $1$ are the only critical values of $\alpha$ and if $\alpha(m)=0$ (resp. $1$) then for any $r>0$, the vector $\partial_r(r,m)$ belongs to $V_2$ (resp $V_1$).
\end{cor}
{\bf Proof.} Let $m\in M$ be  a critical point of $\alpha$. As we already saw, at (\ref{decadix}), it implies that for any $r>0$,
$\hD\hD A^i_{(r,m)}(\partial_r,.)=2\alpha(m)g(\partial_r,.)$. It means that $\partial_r(r,m)$ belongs to a eigenspace of the self-adjoint endomorphism associated to $T_1$ that is  $V_1$ or $V_2$. This, in turn, implies that $\alpha(m)=0$ or $1$.

The manifold $M$ being closed, there exists $(m_+,m_-)\in M^2$  such that $\alpha(m_+)=\max_{m\in M} \alpha(m)$ and  $\alpha(m_-)=\min_{m\in M} \alpha(m)$, therefore $d\alpha(m_{\pm})=0$, $\alpha (m_-)=0$ and $\alpha(m_+)=1$.
$\Box$
\hop
We decompose the vectorfield $\partial_r$ according to the splitting $T\hM=V_1\oplus V_2$. We have 
$$\partial_r=X_1+X_2$$
We decompose the vectors $X_1$ and $X_2$ according to the decomposition $T\hM=\R \partial_r \oplus TM$, we have 
$$X_1=\alpha \partial_r +X, \quad X_2=(1-\alpha)\partial_r - X,$$
 $X$ is therefore a vectorfield on $\hM$ tangent to $M$ satisfying  $\hg(X,X)=\alpha-\alpha^2$. It has a more interesting property:
\begin{prop}[see \cite{Leist}, corollary 4.1]\label{leist}
The vectorfield $2rX$ projects on a vectorfield on $M$ which is the gradient of $\alpha$ (with respect to the metric $g$). 
\end{prop}
{\bf Proof.} Let $(r,m)$ be a point of $\hM$. Let $Z$ be the lift of  vectorfield of $M$ perpendicular at $(r,m)$ to $X$. We known from proposition \ref{GaMa} that $\hD\hD A^1=2T$ therefore, using (\ref{deriv}) at $m$ we have:
$$D\alpha(Z)=\hD\hD A^1(\frac{1}{r}Z,\partial_r)=2\hg(\frac{1}{r}Z,X_1)=0.$$
This means that $Z$ is perpendicular at $m$ to the gradient of $\alpha$ and therefore that $X$ projects on a well-defined direction field of $M$. As $g(2rX,2rX)$ does not depend on $r$ and as $X$ is nowhere lightlike (see corollary \ref{critic}), the vectorfield $2rX$ does project on $M$.

To conclude we just have to show that $D\alpha(2rX)=g(2rX,2rX)$. Using again  (\ref{deriv})  we have:
$$D\alpha(2rX)=\hD\hD A^1(2X,\partial_r)=\hg(2X,X_1)= 2\hg(2X,\partial_r+X)=4\hg(X,X)=g(2rX,2rX).\quad \Box$$
%
\begin{cor}
 The gradient of $A^i$ is the vectorfield $2rX_i$. 
\end{cor}
{\bf Proof.} We have $dA^1=2r\alpha dr+r^2 d\alpha$. Let $v=a\partial_r+h$ be a vector tangent to $\hM$ decomposed according to the splitting $T\hM=\R \partial_r \oplus TM $. We verify that $\hg(2rX_1,.)=dA^1$. Using proposition \ref{leist}, we have 
$$\begin{array}{ll}
dA^1(v)&=2r\alpha a+r^2d\alpha(h)\\
      &= 2r\alpha a+r^2 g(2rX,h)\\
      &= 2r \hg(\alpha \partial_r+X, v).
  \end{array} $$
The situation of $A^1$ and $A^2$ being symmetric the same is also true for $A^2$.
$\Box$ 
\begin{cor}\label{constant}
 The function $A^1$ (resp. $A^2$)  is constant along the the leaves of $\ff_2$ (resp. $\ff_1$).
\end{cor}
{\bf Proof.}  It follows from the fact that $X_i$ belongs to $V_i$ and that $V_1$ and $V_2$ are orthogonal.$\Box$
\begin{prop}\label{moi}
Let $m_0\in U\subset M$ (i.e. such that $\alpha(m)\neq 0,1$), and, for $i\in \{1,2\}$, let $L^i_{(r_0,m_0)}$ be the leaf of $\ff_i$ containing the point $(r,m)$.
Then there exists a point $p$  in $L^1_{(r_0,m_0)}$ (resp.  $p\in L^2_{(r_0,m_0)}$)  such that $\alpha(p)=0$ (resp $\alpha(p)=1$).
\end{prop}
{\bf Proof.}
Let $(r,m)\in L^i_{(r_0,m_0)}$. We are going to look for a critical point of the restriction of  $A^i$ to $L^i_{(r_0,m_0)}$. Classically we follow (backward) the integral curves of the gradient. We note that the gradient of $A^i$ and the gradient of the restriction of $A^i$ are the same. 
Let $(r,\gamma) :\ ]a,b[\ \rightarrow L^i_{(r_0,m_0)}$ be the maximal integral curve of $rX_i$ such that $(r(0),\gamma(0))=(r,m)$.
\begin{lem}\label{minore}
 The image of the  restriction of $\gamma$ to $]a,0]$ lies in a compact set of $\hM$.
\end{lem}
{\bf Proof.} The manifold $M$ being compact, we just have to show that $r(]a,0])$ is contained in a compact subset of $]0,+\infty[$.

We first remark that $rX_1.r=r\alpha\geq 0$ and $rX_2.r=r(1-\alpha)\geq 0$. This implies that $$\forall t\in ]a,0],  r(t)\leq r(0)=r.$$ We thus have a upper bound on $r(t)$.

Moreover, let $(r',m')$ be a point of $L^1_{(r_0,m_0)}$ (resp. of $L^2_{(r_0,m_0)}$) the function $A^{2}$ (resp. $A^1$) being  constant on $L^1_{(r_0,m_0)}$ (resp. on $L^2_{(r_0,m_0)}$), we have $r'^2\alpha(m')=r_0^2\alpha(m_0)\neq 0$ (resp. $r'^2(1-\alpha(m'))=r_0^2(1-\alpha(m_0))\neq 0$). According to corollary \ref{critic}, $0\leq \alpha(m)\leq 1$. Therefore $r'\geq r_0\sqrt {\alpha(m_0)}$ (resp $r'\geq r_0\sqrt {1-\alpha(m_0)}$). As for all $t\in ]a,b[$, we have $(r(t),\gamma(t))\in L^i_{(r_0,m_0)}$, this gives a lower bound for $r(t)$.
$\Box$\hop
It follows from lemma \ref{minore} that there exists a sequence $(t_n)_{n\in\N}$ of points of $]a,0]$ converging to $a$ and such that the sequence  $(\gamma(t_n))_{n\in\N}$ converges in $\hM$ to a point $(r_\infty, m_\infty)$. Let $O$ be  a foliated neighborhood for $\ff_i$  of $(r_\infty, m_\infty)$. There are two possibilities: either $(r_\infty, m_\infty)$ belongs to $L^i_{(r_0,m_0)}$ or the points $\gamma(t_n)$ belong to an infinite number of connected components of $O\cap L^i_{(r_0,m_0)}$ (called plaques). The last case implies that the leaf $L^i_{(r_0,m_0)}$ accumulates around  $(r_\infty, m_\infty)$. As the vector $\partial_r$ is never tangent to $L^i_{(r_0,m_0)}$ this is incompatible with the following straightforward consequence of corollary \ref{constant}.
\begin{fact}
 Let $m\in M$, if the set $\R^*_+\times \{m\}\cap L^i_{(r_0,m_0)}$ contains more than  one point then  $\R^*_+\times \{m\}\subset  L^i_{(r_0,m_0)}$.
\end{fact}
Hence  $(r_\infty, m_\infty)$ is a (local) minimum of the restriction of $A^i$, therefore $X_i(r_\infty, m_\infty)=0$. It implies that $A^i(r_\infty, m_\infty)=0$. If $i=1$ it implies that $\alpha(m_\infty)=0$ and if $i=2$ that $(1-\alpha)(m_\infty)=0$.$\Box$
\hop
Now we use lemma 3.2 of \cite{Ga} or equivalently lemma 6.3 of \cite{Leist}, which says
\begin{prop}
 If a leaf $L^i_{(r,m)}$ contains a point $p$ such that $X_i(p)=0$ then it is flat.
\end{prop}
This proposition together with proposition \ref{moi} entails  that $\hM$ is flat.

It follows from fact \ref{fact} that the curvature of $\hg$ is given by $\widehat R(X,Y)Z=R(X,Y)Z-g(Y,Z)X+g(X,Z)Y$, where $R$ and $\widehat R$ are the curvature of $g$ and $\hg$. It implies that 
$(\hM,\hg)$ is flat if and only if $(M,g)$ has constant curvature equal to $1$.

Moreover, in  \cite{Ma}, Matveev remarked that for any lightlike geodesic $\gamma$ and any  solution $\alpha$ of $(*)$ we have $(\alpha\circ \gamma)'''=0$, therefore solutions of the equation $(*)$ are constant or unbounded along complete lightlike geodesics. If the manifold is closed and complete the function $\alpha$ is constant along lightlike geodesics and therefore constant. It implies that the cone over a closed \emph{complete} pseudo-Riemannian manifold is never decomposable. Therefore we have:
\begin{thm}
 If $(M,g)$ is a closed pseudo-Riemannian (but not Riemannian) manifold with decomposable cone then $(M,g)$ is  \emph{incomplete} and has constant curvature equal to $1$.
\end{thm}
Let us note  that no closed incomplete manifolds of constant curvature are known, probably they do not exist. However the only result in this direction is the  Carri\`ere-Klingler theorem (see  \cite{Ca} and  \cite{Kl}) which asserts that  closed  Lorentzian manifold of constant curvature are complete. Hence we have:
\begin{cor}\label{lorentza}
 If $(M,g)$ is a closed pseudo-Riemannian manifold of signature $(n-1,1)$ or $(1,n-1)$ then its cone is not decomposable.
\end{cor}
In order to extend this result to any signature we need a result weaker than the Carri\`ere-Klingler theorem. 
We need to answer the following.
\begin{quest}
 Does there exist a closed pseudo-Riemannian manifold of signature $(p,q)$ with constant  curvature equal to $1$ whose holonomy group (in the $(G,X)$-structure sense) preserves a non degenerate splitting of $\R^{p+1,q}$ (we recall that the holonomy group lies in the isometry group of the model spaces of constant curvature, ie in $O(p+1,q)$).
\end{quest}
If we assume the manifold complete the answer follows of course from what we have done, but it is not difficult to prove it directly (at least when the pseudo-sphere $\{x\in \R^{p+1,q}\,| \langle x,x\rangle =1\}$ is simply-connected ie when $p>1$ and $q>2$, see \cite{ze} fact 2.3).
\hop
It follows from proposition \ref{para->decomp} that any closed pseudo-Riemannian manifold admitting a non-trivial solution to the equation $(*)$ has a  decomposable cone. Therefore we have proven theorem \ref{pr}. Corollary \ref{lorentza} gives the following
\begin{cor}\label{lorentz}
 If $(M,g)$ is a closed (ie compact without boundary) Lorentzian manifold then any solution of the equation 
$$(*)\qquad DDD \alpha (X,Y,Z)+ 2 (D\alpha  \otimes g)(X,Y,Z)+(D\alpha \otimes g)(Y,X,Z)+(D\alpha \otimes g)(Z,X,Y)=0,$$
is constant.
\end{cor}
\section{Application to the projective Lichnerowicz conjecture.}\label{lichne}
Let $M$ be a  manifold of dimension greater than $2$. Let $g$ and $g'$ be two pseudo-Riemannian metric on $M$. We say that $g$ and $g'$ are \emph{geodesically equivalent} if every $g$-geodesic is a reparametrized $g'$-geodesic. We say that $g$ and $g'$ are \emph{affinely equivalent} if their Levi-Civita connections coincide.

It is interesting to known when geodesic equivalence implies affine equivalence and how big can a space of geodesically equivalent metrics be.  This question is related to the following conjecture (let us note that the conjecture in \cite{KM} assumes completeness instead of compactness).
\hop
{\bf Projective Lichnerowicz conjecture (closed case)} \emph{
Let $G$ be  a connected Lie group acting on a \emph{closed} connected pseudo-Riemannian or Riemannian manifold $(M,g)$ of dimension $n\geq 2$ by projective transformations. Then it acts by affine transformations or $(M,g)$ is a quotient of a Riemannian round sphere.}
\hop
An action is said to be projective if it sends unparametrized geodesics on unparametrized geodesics. The conjecture has been solved in the Riemannian case by V. Matveev in \cite{Ma2} but only  partially solved in the pseudo-Riemannian case (cf the article of V. Kiosak and V. Matveev \cite{KM}).

In \cite{KM} the authors define the \emph{degree of mobility} of a metric $g$ as the dimension of the space of metrics $g'$ which are geodesically equivalent to $g$. This set being linear this number is well defined and never equal to $0$. We will use the following proposition which appears p18 in \cite{KM} as corollary 4.
\begin{prop}[see  \cite{KM}]\label{KM}
 Let $(M,g)$ be a connected pseudo-Riemannian manifold of dimension greater than $2$. If the degree of mobility of $g$ is greater than $2$ and if there exists a metric $g'$ which is geodesically equivalent to $g$ but not affinely equivalent to $g$  then there exists $c\in \R$ such that the equation
$$(**)\qquad DDD \alpha (X,Y,Z)+ c[2 (D\alpha  \otimes g)(X,Y,Z)+(D\alpha \otimes g)(Y,X,Z)+(D\alpha \otimes g)(Z,X,Y)]=0,$$
has a non-constant solution.
\end{prop}
Theorem \ref{pr} allows us to extend the main  theorem of \cite{KM} to the closed case:
\begin{thm}\label{mob}
 Let $(M,g)$ be a \emph{closed} connected pseudo-Riemannian manifold of dimension greater than $2$. If the degree of mobility of $g$ is greater than $2$ and if there exists a metric $g'$ which is geodesically equivalent to $g$ but not affinely equivalent to $g$ then $(M,g)$ is a closed incomplete manifold of constant non zero curvature.
\end{thm}
{\bf Proof.} If $c\neq 0$ we remark that $\alpha$ is a solution of $(**)$ for the pseudo-Riemannian metric $g$ if and only if it is a solution of $(*)$ on the metric $cg$. In that case the theorem follows then from theorem \ref{pr} and proposition \ref{KM}. 

Let us suppose now that  $c=0$. Equation $(**)$ tells us that $\alpha$ is  a function on a closed manifold with parallel Hessian. At a minimum the Hessian must be positive and at a maximum it must be negative. Therefore the Hessian is null, therefore the gradient of $\alpha$ is parallel. But as it vanishes somewhere it vanishes everywhere and $\alpha$ is constant. $\Box$
\hop
The main result of \cite{KM} is almost the same as our theorem \ref{mob} replacing the assumption that the manifold is compact by the one that the metric is complete. It's worth noting that, in that case, the equation $(*)$ does have non-constant solutions. But the very peculiar behavior of solutions of $(*)$ along lightlike geodesic is then  used by the authors to conclude.
We have made one step in the direction of the projective Lichnerowicz conjecture.
\begin{cor}
A closed counter-example to the projective Lichnerowicz conjecture is either a $2$ torus, an incomplete closed manifold of constant non-zero curvature or a metric whose degree of mobility is exactly $2$.
\end{cor}
As above, we can remove the constant curvature case if we suppose that the manifold is Lorentzian.

\hip
\begin{tabular}{ll}
 Address: & Universit\'e Bordeaux 1, Institut de Math\'ematiques de Bordeaux,\\
 &351, cours de la lib\'eration, F-33405 Talence, France\\
E-mail:&{\tt pierre.mounoud@math.u-bordeaux1.fr}
\end{tabular}

\begin{thebibliography}{10}
 \bibitem{Leist} D.V. Alekseevsky, V. Cortes, A.S. Galaev, T. Leistner, \emph{Cones over pseudo-Riemannian manifolds and their holonomy}, to appear in J. Reine Angew. Math. (Crelle's journal). arXiv: 0707.3063v2
\bibitem{Bouba} C. Boubel, \emph{The algebra of the parallel endomorphisms of a pseudo-Riemannian manifold}, in preparation.
\bibitem{Ca} Y. Carri\`ere, \emph{Autour de la conjecture de L. Markus sur les vari\'et\'es affines}, Invent. Math.  95  (1989),  no. 3, 615--628.
\bibitem{CR} Y. Carri\`ere, L. Rozoy, \emph{ Compl\'etude des m\'etriques lorentziennes de $T\sp 2$ et diff\'eormorphismes du cercle.} Bol. Soc. Brasil. Mat. (N.S.)  25  (1994),  no. 2, 223--235.
\bibitem {Ga} S. Gallot, \emph{\'Equations diff\'erentielles caract\'eristiques de la sph\`ere}, Ann. scient. \'Ec. Norm. Sup. 4$^{e}$ s\'erie t.12, 1979, p235-267.
\bibitem{KM} V. Kiosak, V. Matveev, \emph {Proof of projective Lichnerowicz conjecture for pseudo-Riemannian metrics with degree of mobility greater than two}, preprint arXiv:0810.0994v3.
\bibitem{Kl} B. Klingler, \emph{Compl\'etude des vari\'et\'es lorentziennes \`a courbure constante}, Math. Ann.  306  (1996),  no. 2, 353--370.
\bibitem{Ma2} V. Matveev, \emph{Proof of the Lichnerowicz Conjecture}, J. Diff. Geom. 75 (2007).
\bibitem {Ma} V. Matveev, \emph{Gallot-Tanno theorem for pseudo-Riemannian manifolds and a proof that decomposable cones over closed complete pseudo-Riemannian manifolds do not exists.}, preprint arXiv:0906.2410v1
\bibitem{So} A.S. Solodovnikov, \emph{Projective transformations of Riemannian spaces}, Uspehi Mat. Nauk (N.S.) 11(4(70)) (1956) 45--116, MR 0084826, Zbl 0071.15202.
\bibitem{Ta}S. Tanno, \emph{Some differential equations on Riemannian manifolds}, J. Math. Soc. Japan 30 (1978), no. 3, 509--531.
\bibitem{ze} A. Zeghib, \emph{On closed anti de Sitter spacetimes}, Math. Ann. 310(1998), 695-716.
\end{thebibliography}
\end{document}